\documentclass[11pt, epsf]{amsart}
\usepackage{epsfig}
\usepackage{eucal}
\usepackage{mathrsfs}
\usepackage{amssymb}
\usepackage{latexsym}
\usepackage{amsfonts}
\usepackage[latin1]{inputenc}

\newenvironment{namelist}[1]{%
\begin{list}{}
 {
   
   \settowidth{\labelwidth}{#1}
   \setlength{\leftmargin}{1.1\labelwidth}
  }
 }{%
\end{list}}
\newcommand{\bi}{\begin{namelist}}
\newcommand{\ei}{\end{namelist}}

\newcommand{\vt}{\vartriangle}

\newcommand{\p}{\partial}
\newcommand{\no}{\nonumber}

\newtheorem{Def}{Definition}

\newtheorem{Theo}{Theorem}

\newtheorem{Lem}{Lemma}

\title{G-convergence of Dirac operators}
\author{Hasan Almanasreh and Nils Svanstedt}
\thanks{Department of Mathematical Sciences,  University of Gothenburg,
        SE-412 96 Gothenburg, Sweden}
\date{\today}

\begin{document}

\maketitle

\begin{abstract} We consider the linear Dirac operator with a $(-1)$-homogeneous locally periodic potential that varies with respect to a small parameter.
Using the notation of G-convergence for positive self-adjoint operators
in Hilbert spaces we prove G-compactness in the strong resolvent sense for families of projections of Dirac operators. We also prove convergence
of the corresponding point spectrum in the spectral gap.\\\\
%{\bf AMS Subject Classification:} 35B27, 35B40.
\end{abstract}

%-----------------------------------------------------------------------------
\section{Introduction}

In the present work we study the asymptotic behavior of Dirac operators $\tilde{\mathscr{H}}_h$ with respect to a parameter $h\in{\mathbb N}$ as $h\to\infty$.
We consider Dirac operators $\tilde{\mathscr{H}}_h= \tilde{\mathbf{H}}+V_h$ on $L^2(\mathbb {R}^3;\mathbb {C}^4)$, where $\tilde{\mathbf{H}}=\mathbf{H}_0+W+\mathbf{I}$ is a shifted Dirac operator. The operators $\mathbf{H}_0$, $W$, $\mathbf{I}$, and $V_h$ are respectively the free Dirac operator, the Coulomb potential, the $4\!\times\!4$ identity matrix, and a perturbation to $\tilde{\mathbf{H}}$. We will study the asymptotic behavior of $\tilde{\mathscr{H}}_h$ and of the eigenvalues in the gap of the continuous spectrum with respect to
the perturbation parameter $h$.\\

G-convergence theory which deals with convergence of operators,
 is well-known for its applications to homogenization of partial differential equations. The concept
was introduced in the late 60's by De Giorgi and Spagnolo \cite{DEGS, SPA67, SPA75} for linear elliptic and parabolic problems with symmetric coefficients matrices. Later on it was extended to
the non-symmetric case by Murat and Tartar \cite{MUR, TAR1, TAR3} under the name of H-convergence.
A detailed exposition of G-convergence for positive self-adjoint operators is found in Dal Maso \cite{DAL}. In the present
work we will base a lot of our framework on the results in Chapter 12 and 13 in \cite{DAL}. The Dirac operator is unbounded both from above and below. This means that the
theory of G-convergence for positive self-adjoint operators is not directly applicable to Dirac operators. In this work we study self-adjoint projections of Dirac operators which
satisfy the positivity so that the theory of G-convergence becomes applicable.\\

We will consider periodic perturbations, i.e. we will assume that the potential $V_h$ is a periodic function with respect to some regular lattice in $\mathbb{R}^N$. We are then interested in the
asymptotic behaviour of shifted perturbed Dirac operators $\tilde{\mathscr{H}}_h$. This yields homogenization problems for the evolution equation
associated with the Dirac operator $\tilde{\mathscr{H}}_h$
\begin{equation}\no
\left\{ \begin{array}{l}
i\hslash\frac{\partial}{\partial t}\mathbf{u}_h(t,x) =\tilde{\mathscr{H}}_h\mathbf{u}_h(t,x)\,,\\
\mathbf{u}_h(\cdot,0)=\mathbf{u}_h^0\,
\end{array} \right.
\end{equation}
and the corresponding eigenvalue problem
$$
\tilde{\mathscr{H}}_h u_h(x)=\lambda_h u_h(x).
$$

The paper is arranged as follows: In Section 2 we provide the reader with basic preliminaries on Dirac operators, G-convergence and on the concepts needed from spectral theory. In Section 3 we present
and prove the main results.
\section{Preliminaries}
Let $A$ be a linear operator on a Hilbert space. By $\mathbf{R}(A)$, $\mathbf{D}(A)$, and $\mathbf{N}(A)$ we mean the range, domain, and null-space of $A$ respectively.
\subsection{Dirac Operator}
We recall some basic facts regarding the Dirac operator. For more details we refer to the monographs \cite{THA}, \cite{WEI80} and \cite{WEI2003}.

Let $\mathcal{X}$ and $\mathcal{Y}$ denote the Hilbert spaces $H^1(\mathbb {R}^3;\mathbb {C}^4)$ and $L^2(\mathbb {R}^3;\mathbb {C}^4)$, respectively. The free Dirac evolution equation reads
\begin{equation}\label{8}
i\hslash\frac{\partial}{\partial t}\mathbf{u}(t,x) =\mathbf{H}_0\mathbf{u}(t,x)\,,
\end{equation}
where $\mathbf{H}_0:  \mathcal{Y}\longrightarrow \mathcal{Y}$ is the free Dirac operator with domain $\mathbf{D}(H_0)=\mathcal{X}$, which acts on the four-component vector $\mathbf{u}$.
It is a first order linear hyperbolic partial differential equation. The free Dirac operator $\mathbf{H}_0$ has the form
\begin{equation}\label{9}
\mathbf{H}_0 = -i{\hslash}c\boldsymbol{\alpha}\cdot \mathbf{\nabla} + mc^2 \beta\,.
\end{equation}
Here $\boldsymbol{\alpha}\cdot \mathbf{\nabla} = \sum_{i=1}^3 \alpha_i \frac{\p}{\p x_i}$,
$\hslash$ is the Planck constant divided by $2\pi$, the constant $c$ is the speed of light, $m$ is the particle rest mass and $\boldsymbol{\alpha}  = (\alpha_1,\alpha_2,\alpha_3)$ and $\beta$ are the $4\!\times\!4$ Dirac matrices given by
$$
\alpha_i = \left(
\begin{array}{cc}
0 & \sigma_i \\
\sigma_i & 0
\end{array}
\right)\;\;\text{and}\;\;
\beta = \left(
\begin{array}{cc}
I & 0 \\
0 & -I
\end{array}
\right)\, .
$$
Here $I$ and $0$ are the $2\!\times\!2$ unity and zero matrices, respectively, and the $\sigma_i$'s are the $2\!\times\!2$ Pauli matrices
$$
\sigma_1 = \left(
\begin{array}{cc}
0 & 1 \\
1 & 0
\end{array}
\right),\;\;
\sigma_2 = \left(
\begin{array}{cc}
0 & -i \\
i & 0
\end{array}
\right)
\;\;\text{and}\;\;
\sigma_3 = \left(
\begin{array}{cc}
1 & 0 \\
0 & -1
\end{array}
\right)\, .
$$
Note that a separation of variables in (\ref{8}) yields the Dirac eigenvalue problem
\begin{equation}\label{10}
\mathbf{H}_0u(x)=\lambda u(x)\, ,
\end{equation}
where $u(x)$ is the spatial part of the wave function $\mathbf{u}(x,t)$ and $\lambda$ is the total energy of the particle. The free Dirac operator $\mathbf{H}_0$ is essentially self-adjoint on
$C^\infty_0(\mathbb{R}^3;\mathbb{C}^4)$ and self-adjoint on $\mathcal{X}$. Moreover, its spectrum, $\sigma(\mathbf{H}_0)$, is purely absolutely continuous (i.e. its spectral measure is absolutely continuous with respect to the Lebesgue measure) and given by
$$
\sigma(\mathbf{H}_0) = (-\infty,-mc^2]\cup[mc^2,+\infty)\, .
$$
$\mathbf{H}_0$ describes the motion of an electron that moves freely without external force. Let us now introduce an external field given by a $4\!\times\!4$ matrix-valued function $W$,
$$
W(x)=W_{ij}(x) \;\;\;\; i,j=1,2,3,4.
$$
It acts as a multiplication operator in $L^2(\mathbb{R}^3;\mathbb{C}^4)$, thus the free Dirac operator with additional field $W$ is of the form
\begin{equation}\label{11}
\mathbf{H}=\mathbf{H}_0+W\, .
\end{equation}
The operator $\mathbf{H}$ is essentially self-adjoint on $C^\infty_0(\mathbb{R}^3;\mathbb{C}^4)$ and self-adjoint on the Sobolev space $\mathcal{X}$ provided that $W$ is Hermitian and satisfies the following estimate (see e.g. \cite{THA})
\begin{equation}\label{5}
|W_{ij}(x)|\leq a\frac{c}{2|x|}+b\,,\quad \forall x\in \mathbb{R}^3\backslash\{0\}\quad i,j=1,2,3,4,
\end{equation}
the constant $c$ is the speed of light, $a<1$, and $b>0$. From now on we let $W(x)$ be the Coulomb potential $W(x)=\frac{-Z}{x}\mathbf{I}$, where $Z$ is the elctric charge number (without ambiguity, $\mathbf{I}$ is usually dropped from the Coulomb term for simplicity). The spectrum of the Dirac operator with Coulomb potential
is given by
$$
\sigma(\mathbf{H}) = (-\infty,-mc^2]\cup\{\lambda^k\}_{k\in\mathbb N}\cup[mc^2,+\infty),
$$
where $\{\lambda^k\}_{k\in\mathbb N}$ is a discrete sequence of eigenvalues in the "gap" and the remaining part of the
spectrum is the continuous part $\sigma(\mathbf{H}_0)$.

In the present paper we consider a parameter-dependent perturbation added to the Dirac operator with Coulomb potential. The purpose is to investigate the asymptotic behavior of the corresponding eigenvalues in the gap and the convergence properties. To this end we introduce a $4\!\times\!4$ matrix-valued function $V_h = V_h(x)$ and define the operator $\mathscr{H}_h$ as
\begin{equation}\label{12}
\mathscr{H}_h=\mathbf{H}+V_h\, .
\end{equation}
We recall that a function $F$ is called homogeneous of degree $p$ if for any nonzero scalar $a$, $F(ax)=a^pF(x)$. The next theorem is of profound importance
for the present work.
\begin{Theo}
Let $W$ be Hermitian and satisfy the bound (\ref{5}) above. Further, for any fixed $h\in\mathbb{N}$, let $V_h$ be a measurable $(-1)$-homogeneous Hermitian $4\!\times\!4$ matrix-valued function
with entries in $L^p_{loc}(\mathbb{R}^3)$, $p>3$. Then $\mathscr{H}_h$ is essentially self-adjoint on $C^\infty_0(\mathbb{R}^3;\mathbb{C}^4)$ and self-adjoint on $\mathcal{X}$. Moreover
$$
\sigma(\mathscr{H}_h) = (-\infty,-mc^2]\cup\{\lambda_h^k\}_{k\in\mathbb N}\cup[mc^2,+\infty),
$$
where $\{\lambda_h^k\}_{k\in\mathbb N}$ is a discrete sequence of parameter dependent eigenvalues corresponding to the Dirac eigenvalue problem $\mathscr{H}_h u_h(x)=\lambda_h u_h(x)$.
\end{Theo}
\underline{Proof}. \emph{See \cite{WEI80,WEI2003}}.\hfill{$\blacksquare$}
\\\\
We will as a motivating example consider perturbations which are locally periodic and of the form
$V_h(x) = V_1(x)V_2(hx)$. The entries of $V_1$ are assumed to be $(-1)$-homogeneous. The entries of $V_2(y)$ are assumed to be periodic with respect to a regular
lattice in $\mathbb{R}^3$. This can also be phrased that they are defined on the unit torus $\mathbb{T}^3$.\\\\
The evolution equation associated with the Dirac operator $\mathscr{H}_h$ reads
\begin{equation}\label{14}
\left\{ \begin{array}{l}
i\hslash\frac{\partial}{\partial t}\mathbf{u}_h(t,x) =\mathscr{H}_h\mathbf{u}_h(t,x)\,,\\
\mathbf{u}_h(\cdot,0)=\mathbf{u}_h^0\,.
\end{array} \right.
\end{equation}
By the Stone theorem, since $\mathscr{H}_h$ is self-adjoint on $\mathcal{X}$, there exists a unique solution
$\mathbf{u}_h$ to (\ref{14}) given by
\begin{equation}\label{15}
\mathbf{u}_h(\cdot,t)=\mathcal{U}_h(t)\mathbf{u}_h^0\, , \;\forall \mathbf{u}_h^0\in \mathcal{X}\,,
\end{equation}
where $\mathcal{U}_h(t)=\exp(-(i/\hslash)\mathscr{H}_ht)$ is the strongly continuous unitary operator generated by the infinitesimal operator $-(i/\hslash)\mathscr{H}_h$ on $\mathcal{Y}$, see e.g. \cite{KAT} or \cite{THA}.

In the sequel we will consider a shifted family of Dirac operators denoted by $\tilde{\mathscr{H}}_h$ and defined as $\tilde{\mathscr{H}}_h=\tilde{\mathbf{H}}+V_h$, where $\tilde{\mathbf{H}}=\mathbf{H}+mc^2\mathbf{I}$. Also without loss of generality we will in the sequel put $\hslash=c=m=1$. By Theorem 1, for any $h\in\mathbb{N}$, we then get
$$
\sigma(\tilde{\mathscr{H}}_h) = (-\infty,0]\cup\{\tilde{\lambda}_h^k\}_{k\in\mathbb N}\cup[2,\infty).
$$
\subsection{G-convergence}
For more detailed information on G-convergence we refer to e.g. \cite{DEF, SVA99} for the application to elliptic and parabolic
partial differential operators, and to the monograph \cite{DAL} for the application to general self-adjoint operators. Here we recall some basic facts about
G-convergence for self-adjoint operators in $\mathcal{Y}$.

In the present work we frequently write $A_h$ converges to $A$ when we mean that the sequence $\{A_h\}$ converges to $A$. Let $\lambda\geq0$, by $\mathcal{P}_\lambda (\mathcal{Y})$ we denote the class of self-adjoint operators $A$ on a closed linear
subspace $\mathscr{V}=\overline{\mathbf{D}(A)}$ of $\mathcal{Y}$ such that $\langle Au,u\rangle\geq\lambda ||u||^2_\mathcal{Y}\; \forall u\in \mathbf{D}(A)$.
\begin{Def}
Let $\lambda>0$, and let $\{A_h\}\subset\mathcal{P}_\lambda(\mathcal{Y})$ then we say that $A_h$ G-converges to $A\in\mathcal{P}_\lambda(\mathcal{Y})$, denoted
$A_h\xrightarrow[\text{\tiny G,w}]{\;\,\text{\tiny G,s}\;\,}A$ in $\mathcal{Y}$ if $A_h^{-1}P_hu\xrightarrow[\text{\tiny w}]{\;\,\text{\tiny s}\;\,}A^{-1}P u$ in $\mathcal{Y}$ $\forall u\in \mathcal{Y}$, where $s$ and $w$ refer to strong and weak topologies respectively, and $P_h$ and $P$ are the orthogonal projections onto $\mathscr{V}_h:=\overline{\mathbf{D}(A_h)}$ and  $\mathscr{V}:=\overline{\mathbf{D}(A)}$ respectively. Also we say $\{A_h\}\subset\mathcal{P}_0(\mathcal{Y})$ converges to $\{A\}\subset\mathcal{P}_0(\mathcal{Y})$ in the strong resolvent sense $(SRS)$ if $(\lambda I+A_h)\xrightarrow[]{\;\,\text{\tiny G,s}\;\,}(\lambda I+A)$ in $\mathcal{Y}$ $\forall\lambda>0$.
\end{Def}

The following result provides a useful criterion for G-convergence of self-adjoint operators. See \cite{DAL} for a proof.
\begin{Lem}
Given $\lambda>0$, $\{A_h\}\subset\mathcal{P}_\lambda(\mathcal{Y})$ and the orthogonal projection $P_h$ onto $\mathscr{V}_h$. Suppose
that for every $u \in\mathcal{Y}$,
$A_h^{-1}P_hu$ converges strongly (resp. weakly) in $\mathcal{Y}$, then
there exists an operator $A\in\mathcal{P}_\lambda(\mathcal{Y})$ such that $A_h\xrightarrow{\;\,\text{\tiny G,s}\;\,}A$ (resp. $A_h\xrightarrow{\;\,\text{\tiny G,w}\;\,}A$) in $\mathcal{P}_\lambda(\mathcal{Y})$.
\end{Lem}

From now on we will just use the word "converge" instead of saying "strongly converge", hence $A_h\xrightarrow{\;\,\text{\tiny G}\;\,}A$ instead of $A_h\xrightarrow{\;\,\text{\tiny G,s}\;\,}A$.

\subsection{Some Basic Results in Spectral Theory}
For more details see \cite{BIR,KAT,WEI80}. Given a Hilbert space $X$, let $(\mathscr{U},\mathscr{A})$ be a measurable space for $\mathscr{U}\subseteq\mathbb{C}$ and $\mathscr{A}$ being a $\sigma-$algebra on $\mathscr{U}$. Let $\mathbb{P}^{X}=\mathbb{P}(X)$ be the set of orthogonal projections on $X$, then $E:\mathscr{A}\longrightarrow\mathbb{P}^{X}$ is called a spectral measure if it satisfies the following
\begin{itemize}
\item [$(i)$] $E(\emptyset)=\mathbf{0}$ (This condition is superfluous given the next properties).
\item [$(ii)$] Completeness; $E(\mathscr{U})=\mathbf{I}$.
\item [$(iii)$] Countable additivity; if $\{\vartriangle_n\}\subset\mathscr{A}$ is a finite or a countable set of disjoint elements and $\vartriangle=\displaystyle\cup_n\!\vt_n$, then $E(\vt)=\displaystyle\sum_n E(\vt_n)$.
\end{itemize}

If $E$ is spectral measure then $E(\vt_1\cap\vt_2)=E(\vt_1)E(\vt_2)=E(\vt_2) E(\vt_1)$, also $E$ is modular, i.e. $E(\vt_1\cup\vt_2)+E(\vt_1\cap\vt_2)=E(\vt_1)+ E(\vt_2)$. For an increasing sequence of sets $\vt_n$, $\displaystyle \lim_{n\to\infty}E(\vt_n)=E(\displaystyle\cup_n\!\!\vt_n)$, while if $\vt_n$ is a decreasing sequence then $\displaystyle \lim_{n\to\infty}E(\vt_n)=E(\displaystyle\cap_n\!\vt_n)$. Because of the idempotence property of the spectral measure we have $\|E_n u\|^2=\langle E_n u,E_n u \rangle=\langle E_n^2 u,u \rangle=\langle E_n u,u \rangle\to\langle Eu,u \rangle=\|Eu\|^2$, which means that the weak convergence and the strong convergence of a sequence of spectral measures $E_n$ are equivalent.

Let $E_{u,u}(\vt)$ be the finite scalar measure on $\mathscr{A}$ generated by $E$,
\begin{equation*}
E_{u,u}(\vt)=\langle E(\vt)u,u \rangle=\|E(\vt)u\|^2
\end{equation*}
and $E_{u,v}(\vt)$ be the complex measure
\begin{equation*}
E_{u,v}(\vt)=\langle E(\vt)u,v \rangle\,,\; \forall u,v\in X\,.
\end{equation*}
By the above notations $E_{u,v}(\vt)\leq \|E(\vt)u\|\|E(\vt)v\| \leq \|u\|\|v\|$.

Let $\mathscr{U}=\mathbb{R}$. The spectral measure on the real line corresponding to an operator $S$ is denoted by $E^S(\lambda)$ (where the superscript $S$ indicates that the spectral measure $E$ corresponds to a specific operator $S$)
\begin{equation*}
E^S(\lambda)=E^S(\vt)\,,\; \text{where }\vt=(-\infty,\lambda)\,,\;\text{for }\lambda\in\mathbb{R}\,.
\end{equation*}
Clearly $E^S(\lambda)$ is monotonic (nondecreasing), i.e. $E^S(\lambda_1)\leq E^S(\lambda_2)$ for $\lambda_1\leq\lambda_2$. Also $\displaystyle\lim_{\lambda\to-\infty}E^S(\lambda)=\mathbf{0}$ and $\displaystyle\lim_{\lambda\to\infty}E^S(\lambda)=\mathbf{I}$. $E^S(\lambda)$ is self-adjoint, idempotent, positive, bounded, right continuous operator ($\displaystyle\lim_{t\to0^+}E^S(\lambda+t)=E^S(\lambda)$), and discontinuous at each eigenvalue of the spectrum. If $\lambda$ is an eigenvalue, then we define $p(\lambda)=E^S(\lambda)-E^S(\lambda-0)$ to be the point projection onto the eigenspace of $\lambda$. For $\lambda$ being in the continuous spectrum $p(\lambda)=\mathbf{0}$.

Now we state the spectral theorem for self-adjoint operators.
\begin{Theo}
For a self-adjoint operator $S$ defined on a Hilbert space $X$ there exists a unique spectral measure $E^S$ on $X$ such that
\begin{itemize}
\item [$(i)$] $S=\displaystyle\int_{\sigma(S)}\lambda\,dE^S(\lambda)$.
\item [$(ii)$] $E(\vt)=\mathbf{0}$ if $\vt\!\cap\,\sigma(S)=\emptyset$.
\item [$(iii)$] If $\vt\subset \mathbb{R}$ is open and $\vt\!\cap\sigma(S)\neq \emptyset$, then $E(\vt)\neq \mathbf{0}$.
\end{itemize}
\underline{Proof}. \emph{See e.g. \cite{BIR}}.\hfill{$\blacksquare$}
\end{Theo}

\section{The main results}

Consider the family $\{\tilde{\mathscr{H}}_h\}_{h\in\mathbb{N}}$ of Dirac operators with domain $\mathbf{D}(\tilde{\mathscr{H}}_h)=\mathcal{X}$. We will state and prove some useful theorems for operators of the class $\mathcal{P}_\lambda(\mathcal{Y})$ for $\lambda\geq 0$, where $\mathcal{X}$ and $\mathcal{Y}$ are the Hilbert spaces defined above.
The theorems are valid for general Hilbert spaces.\\

The following theorem gives a bound for the inverse of operators of the class $\mathcal{P}_\lambda(\mathcal{Y})$ for $\lambda > 0$.
\begin{Theo}
Let $A$ be a positive and self-adjoint operator on $\mathcal{Y}$ and put $B=A+\lambda I$. Then for $\lambda >0$
\begin{itemize}
\item [($i$)] $B$ is injective. Moreover, for every $v\in\mathbf{R}(B)$, $\langle B^{-1}v,v\rangle\geq \lambda \|B^{-1}v\|^2_\mathcal{Y}$ and $\|B^{-1}v\|_\mathcal{Y}\leq \lambda^{-1}\|v\|_\mathcal{Y}$.
\item [($ii$)] $\mathbf{R}(B)=\mathcal{Y}$.
\end{itemize}
\end{Theo}
\hspace{-4mm}\underline{\emph{Proof}}. See Propositions 12.1 and 12.3 in \cite{DAL}.
\\\\
The connection between the eigenvalue problems of the operator and its G-limit of the class $\mathcal{P}_\lambda (\mathcal{Y})$ for $\lambda\geq0$ is addressed in the next two theorems. Here we prove the critical case when $\lambda=0$, where for $\lambda>0$ the proof is analogous and even simpler.
\begin{Theo}
Given a family of operators $\{A_h\}$ of the class $\mathcal{P}_0(\mathcal{Y})$ G-converging to $A\in\mathcal{P}_0(\mathcal{Y})$ in the strong resolvent sense. Let $u_h$ be the solution of $A_hu_h=f_h$, where $\{f_h\}$ is converging to $f$ in $\mathcal{Y}$. If $\{u_h\}$ converges to $u$ in $\mathcal{Y}$, then $u$ solves the G-limit problem $Au=f$.
\end{Theo}
\hspace{-4mm}\underline{\emph{Proof}}. Since $A_h$ G-converges to $A$ in the strong resolvent sense
\begin{equation}\label{59}
B_h^{-1}P_h v\rightarrow B^{-1}P v\,,\; \forall v\in \mathcal{Y}\,,
\end{equation}
where $B_h$ and $B$ are $A_h+\lambda I$ and $A+\lambda I$ respectively. Note that by Theorem 4, $\mathbf{D}(B_h^{-1})=\mathbf{R}(B_h)=\mathcal{Y}$, so the projections $P_h$ and $P$ are unnecessary. \\
Consider $A_hu_h=f_h$ which is equivalent to $B_hu_h=f_h+\lambda u_h$, by the definition of $B_h$ we have $u_h=B_h^{-1}(f_h+\lambda u_h)$. Define $\mathscr{J}_h=f_h+\lambda u_h$ which is clearly convergent to $\mathscr{J}=f+\lambda u$ in $\mathcal{Y}$ by the assumptions. Therefore $B_h^{-1} \mathscr{J}_h\rightarrow B^{-1} \mathscr{J}$, this is because
\begin{eqnarray*}
\begin{array}{ll}
\|B_h^{-1}\mathscr{J}_h-B^{-1}\mathscr{J}\|_\mathcal{Y}&\!\!\!=\|B_h^{-1}\mathscr{J}_h-B_h^{-1} \mathscr{J}+B_h^{-1}\mathscr{J} -B^{-1}\mathscr{J}\|_\mathcal{Y}\\
&\!\!\!\leq\|B_h^{-1}\|_\mathcal{Y}\,\|\mathscr{J}_h-\mathscr{J}\|_\mathcal{Y}+\|B_h^{-1} \mathscr{J}-B^{-1}\mathscr{J}\|_\mathcal{Y}\\
&\;\;\longrightarrow 0.
\end{array}
\end{eqnarray*}
The convergence to zero follows with help of (\ref{59}) and the boundedness of the inverse operator $B_h^{-1}$. Thus, for all $v\in\mathcal{Y}$
\begin{equation*}
\displaystyle\langle u,v\rangle=\displaystyle\lim_{h\to\infty}\langle u_h,v\rangle=\displaystyle\lim_{h\to\infty}\langle B_h^{-1}\mathscr{J}_h,v\rangle=
\langle B^{-1}\mathscr{J},v\rangle\,.
\end{equation*}
Hence $\langle u-B^{-1}\mathscr{J},v\rangle=0$ for every $v\in \mathcal{Y}$, which implies $Bu=\mathscr{J}$, therefore $Au=f$.\hfill{$\blacksquare$}
\begin{Theo}
Let $\{A_h\}$ be a sequence in $\mathcal{P}_0(\mathcal{Y})$ which G-converges to $A\in\mathcal{P}_0(\mathcal{Y})$ in the strong resolvent sense,
and let $\{\mu_h,u_h\}$ be the solution of the eigenvalue problem $A_hu_h=\mu_hu_h$. If
$\{\mu_h,u_h\}\rightarrow \{\mu,u\}$ in $\mathbb{R}\times\mathcal{Y}$, then the limit couple $\{\mu,u\}$ is the solution of the eigenvalue problem $Au=\mu u$.
\end{Theo}
\hspace{-4mm}\underline{\emph{Proof}}. The proof is straight forward by assuming $f_h=\mu_hu_h$ (which converges to $\mu u$ in $\mathcal{Y}$) in the previous theorem.\hfill{$\blacksquare$}\\

The convergence properties of self-adjoint operators has quite different implications on the asymptotic behavior of the spectrum, in particular on the asymptotic behavior
of the eigenvalues, depending on the type of convergence. For a  sequence $\{A_h\}$  of operators which converges uniformly to a limit operator $A$ nice results can be drawn for the spectrum. Exactly speaking $\{\sigma(A_h)\}$ converges to $\sigma(A)$ including the isolated eigenvalues. The same conclusion holds if the uniform convergence is replaced by the uniform resolvent convergence, see e.g. \cite{KAT}. In the case of strong convergence (the same for strong resolvent convergence), if the sequence $\{A_h\}$ is strongly convergent to $A$, then every $\lambda\in\sigma(A)$ is the limit of a sequence $\{\lambda_h\}$ where $\lambda_h \in\sigma(A_h)$, but not the limit of every such sequence $\{\lambda_h\}$ lies in the spectrum of $A$, (see the below example taken from \cite{WEI97}). For weakly convergent sequences of operators no spectral implications can be extracted.\\\\
\textbf{Example.} Let $A_{i,h}$ be an operator in $L^2(\mathbb{R})$ defined by
\begin{equation*}
A_{i,h}=-\frac{d^2}{dx^2}+V_{i,h}(x)\,,\quad \text{for }h\in\mathbb{N}\text{ and }i=1,2\,,
\end{equation*}
where
\begin{eqnarray*}
\begin{array}{ll}
V_{1,h}(x)=\left\{\begin{array}{ll}
-1\,,&\text{if }h\leq x\leq h+1\,,\\
0\,,&\text{Otherwise}\,,
\end{array}\right.&\text{and}\quad
V_{2,h}(x)=\left\{\begin{array}{ll}
-1\,,&\text{if }x\geq h\,,\\
0\,,&\text{Otherwise}\,.
\end{array}\right.
\end{array}
\end{eqnarray*}
The operator $A_{i,h}$ converges to $A=-\frac{d^2}{dx^2}$ in the strong resolvent sense as $\,h\to \infty$ for both $i=1,2$. One can compute the spectrum for the three operators and obtain $\sigma(A_{1,h}) =[0,\infty)\cup\{\mu_h\}$ for $\mu_h$ being a simple eigenvalue in $[-1,0]$ and $\sigma(A_{2,h}) =[-1,\infty)$, whereas for the unperturbed limit operator $A$ the spectrum consists of just the continuous spectrum, i.e. $\sigma(A)=[0,\infty)$.\\

Since the uniform convergence is not always the case for operators, the theorem below provides some criteria by which the G-convergence of an operator in the set $\mathcal{P}_\lambda(\mathcal{Y})$ (and hence the G-convergence in the strong resolvent sense of operators of the class $\mathcal{P}_0(\mathcal{Y})$) implies the convergence of the corresponding eigenvalues.
\begin{Theo}
Let $\{A_h\}$ be a family of operators in $\mathcal{P}_\lambda(\mathcal{Y})$,
$\lambda>0$, with domain $\mathcal{X}$. If $A_h$ G-converges to $A\in\mathcal{P}_\lambda(\mathcal{Y})$, then $\mathcal{K}_h:=A_h^{-1}$ converges in the norm of $\mathcal{B}(\mathcal{Y})$ ($\mathcal{B}(\mathcal{Y})$ is the set of bounded linear operators on $\mathcal{Y}$) to $\mathcal{K}:=A^{-1}$. Moreover the $k^{th}$ eigenvalue $\mu_h^k$ of $A_h$ converges to the $k^{th}$ eigenvalue $\mu^k$ of $A$
and the associated $k^{th}$ eigenvector $u_h^k$ converges to $u^k$ weakly in $\mathcal{X}$, $\forall k\in \mathbb{N}$.
\end{Theo}
\hspace{-4mm}\underline{\emph{Proof}}. By the definition of supremum norm
\begin{equation}\label{62}
\|\mathcal{K}_h-\mathcal{K}\|_{\mathcal{B}(\mathcal{Y})}= \displaystyle \sup_{\|v\|_\mathcal{Y}=1}\| \mathcal{K}_hv-\mathcal{K}v\|_{\mathcal{Y}}=  \displaystyle \sup_{\|v\|_\mathcal{Y}\leq1}\| \mathcal{K}_hv-\mathcal{K}v\|_{\mathcal{Y}}\,.
\end{equation}
Also, by the definition of supremum norm there exists a sequence $\{v_h\}\subset\mathcal{Y}$ with $\|v_h\|_\mathcal{Y}\leq1$ such that
\begin{equation}\label{63}
\|\mathcal{K}_h-\mathcal{K}\|_{\mathcal{B}(\mathcal{Y})}\leq\| \mathcal{K}_hv_h-\mathcal{K}v_h\|_{\mathcal{Y}} +\frac{1}{h}\,.
\end{equation}
It is well-known that $\mathcal{K}_h$ and $\mathcal{K}$ are compact self-adjoint operators on $\mathcal{Y}$. Both are bounded operators,
by Theorem 3, with compact range $\mathcal{X}$ of $\mathcal{Y}$.\\
Consider now the right hand side of (\ref{63}). We write this as
\begin{equation} \no
\| \mathcal{K}_hv_h-\mathcal{K}v_h\|_{\mathcal{Y}} +\frac{1}{h} \leq \|\mathcal{K}_hv_h-\mathcal{K}_hv\|_{\mathcal{Y}} +\|\mathcal{K}_hv-\mathcal{K}v\|_{\mathcal{Y}}+
\|\mathcal{K}v_h-\mathcal{K}v\|_{\mathcal{Y}}+\frac{1}{h}.
\end{equation}
The first and the third terms converge to zero by the compactness of $\mathcal{K}_h$ and $\mathcal{K}$ on $\mathcal{Y}$  and the second term converges to zero by the G-convergence of $A_h$ to $A$.
Consequently
\begin{equation}\label{K}
\|\mathcal{K}_h-\mathcal{K}\|_{\mathcal{B}(\mathcal{Y})} \to 0.
\end{equation}
Consider the eigenvalue problems associated to $A_h^{-1}$ and $A^{-1}$
\begin{equation}\label{68}
A_h^{-1}v_h^k=\lambda_h^kv_h^k\;,\; k\in \mathbb{N}
\end{equation}
and
\begin{equation}\label{69}
A^{-1}v^k=\lambda^kv^k\;,\; k\in \mathbb{N}\,.
\end{equation}
Since $A_h^{-1}$ and $A^{-1}$ are compact and self-adjoint operators it is well-known that there exists
infinite sequences of eigenvalues $\lambda_h^1\geq\lambda_h^2\geq\cdots$ and $\lambda^1\geq\lambda^2\geq\cdots$ accumulating at the origin, respectively.
Define $\mu_h^k:=(\lambda_h^k)^{-1}$ and $\mu^k:=(\lambda^k)^{-1}$ for all $k\in \mathbb{N}$.
Consider now the spectral problems associated to $A_h$ and $A$
\begin{equation}\label{66}
A_hu_h^k=\mu_h^ku_h^k\;,\; k\in \mathbb{N}
\end{equation}
and
\begin{equation}\label{67}
Au^k=\mu^ku^k\;,\; k\in \mathbb{N}\,.
\end{equation}
There exists infinite sequences of eigenvalues $0<\mu_h^1\leq\mu_h^2\leq\cdots$ and $0<\mu^1\leq\mu^2\leq\cdots$ respectively. By the compactness of $\mathcal{K}_h$ and $\mathcal{K}$ the sets
$\{\lambda_h^k\}_{k=1}^\infty$ and $\{\lambda^k\}_{k=1}^\infty$ are bounded in $\mathbb{R}$, thus the proof is complete by virtue of the following lemma.\hfill{$\blacksquare$}
\begin{Lem}
Let $\mathcal{X}$, $\mathcal{Y}$, $\mathcal{K}_h$, $\mathcal{K}$, $\lambda_h^k$ and $\lambda^k$ be as in Theorem 6, and let $A_h\in \mathcal{P}_\lambda(\mathcal{Y})$, $\lambda>0$. There is a sequence $\{r_h^k\}$ converging to zero with $0<r_h^k<\lambda^k$ such that
\begin{equation}\label{70}
|\lambda_h^k-\lambda^k|\leq c\frac{\lambda^k}{\lambda^k-r_h^k} \displaystyle\sup_{\substack{u\in\mathscr{N}(\lambda^k,\mathcal{K})\\\|u\|_{\mathcal{Y}}=1}} \|\mathcal{K}_hu-\mathcal{K}u\|_ \mathcal{Y}\,,
\end{equation}
where $c$ is a constant independent of $h$, and $\mathscr{N}(\lambda^k,\mathcal{K})= \{u\in\mathbf{D}(\mathcal{K})\,;\; \mathcal{K}u=\lambda^k u\}$ is the eigenspace of $\mathcal{K}$ corresponding to $\lambda^k$.
\\\\
\underline{Proof}. See Theorem 1.4 and Lemma 1.6 in \cite{OLE} Chapter 3.\hfill{$\blacksquare$}
\end{Lem}
We can now complete the proof of Theorem 6. By the G-convergence of $A_h$ to $A$ we obtain, by using Lemma 2 and (\ref{K}), convergence of the eigenvalues and eigenvectors, i.e. $\mu_h^k \to \mu^k$ and $u_h^k \to u^k$ weakly in $\mathcal{X}$ as $h\to \infty$.\hfill{$\blacksquare$}\\

Let us now return to the shifted and perturbed Dirac operator $\tilde{\mathscr{H}}_h$. We will throughout this section assume the hypotheses of Theorem 1. We further assume that the
$4\times 4$ matrix-valued function $V_h$ is of the form $V_h(x) = V_1(x)V_2(hx)$ where $V_1$ is (-1)-homogeneous and where the entries of $V_2(y)$ are 1-periodic in $y$, i.e.
$$
V_2^{ij}(y+k) = V_2^{ij}(y),\;\; k\in \mathbb{Z}^3.
$$
We also assume that the entries of $V_2$ belong to $L^\infty(\mathbb{R}^3)$. It is then well-known that
\begin{equation}\label{meanvalue}
V_2^{ij}(hx) \to M(V_2^{ij}) = \int_{\mathbb{T}^3} V_2^{ij}(y)\, dy,
\end{equation}
in $L^\infty(\mathbb{R}^3)$ weakly{*}, where $\mathbb{T}^3$ is the unit torus in $\mathbb{R}^3$. It easily also follows from this mean-value property that
$$
V_h \to V_1M(V_2),
$$
in $L^p(\mathbb{R}^3)$ weakly for $p > 3$, cf the hypotheses in Theorem 1.\\\\
We are now interested in the asymptotic behavior of the operator and the spectrum of the perturbed Dirac operator $\tilde{\mathscr{H}}_h$.
We recall the spectral problem for $\tilde{\mathscr{H}}_h$, i.e.
$$
\mathscr{\tilde H}_h u_h(x) = \tilde{\lambda}_h u_h(x)
$$
where there exists a discrete set of eigenvalues $\{\tilde{\lambda}_h^k\}$, $k=1,\, 2,\ldots$ and a corresponding set of mutually orthogonal eigenfunctions $\{u_h^k\}$.
We know, by Theorem 1, that the eigenvalues (or point spectrum)  $\sigma^p(\tilde{\mathscr{H}}_h) \subset (0,2)$. We also know that $\tilde{\mathscr{H}}_h$ has a
continuous spectrum $\sigma^c(\tilde{\mathscr{H}}_h)=(-\infty,0]\cup[2,\infty)$. This means that the Dirac operator is neither a positive or negative (semi-definite) operator
and thus the G-convergence method introduced in the previous section for positive self-adjoint operators is not directly applicable. In order to use G-convergence
methods for the asymptotic analysis of $\tilde{\mathscr{H}}_h$ we therefore use spectral projection and study the corresponding asymptotic behavior of projections
$\tilde{\mathscr{H}}_h$ which are positive so that G-convergence methods apply.\\\\
Let $\mathscr{A}$ be a fixed $\sigma$-algebra of subsets of $\mathbb{R}$, and let $(\mathbb{R},\mathscr{A})$ be a measurable space.
Consider the spectral measures $E^{\tilde{\mathbf{H}}}$ and $E^{\tilde{\mathscr{H}}_h}$ of the families of Dirac operators $\tilde{\mathscr{H}}_h$ and $\tilde{\mathbf{H}}$ respectively, each one of these measures maps $\mathscr{A}$ onto $\mathbb{P}^\mathcal{X}$, where $\mathbb{P}^\mathcal{X}$ is the set of orthogonal projections on $\mathcal{X}$. By the spectral theorem
\begin{equation}\label{71half}
\tilde{\mathscr{H}}_h = \displaystyle\int_{\sigma(\tilde{\mathscr{H}}_h)}\lambda\,dE^{\tilde{\mathscr{H}}_h}(\lambda).
\end{equation}
By the spectral theorem we can also write
\begin{equation}\label{71half2}
\displaystyle\int_{\sigma(\tilde{\mathbf{H}})}\lambda\,dE^{\tilde{\mathbf{H}}}(\lambda) + V_h,
\end{equation}
since $V_h$ is a multiplication operator.\\\\
We recall that $\mathbf{D}(\tilde{\mathscr{H}}_h) = \mathcal{X}$, let now
$$
\mathscr{N}_h^k= \{u_h\in \mathcal{X};\tilde{\mathscr{H}}_h u_h=\lambda^k_h u_h\},
$$
i.e. the eigenspace of $\tilde{\mathscr{H}}_h$ corresponding to the eigenvalue $\lambda_h^k$.
Further define the sum of mutual orthogonal eigenspaces
$$
\mathcal{X}_h^p = \oplus_{k\in\mathbb{N}}\mathscr{N}_h^k,
$$
where $\mathcal{X}_h^p$ is a closed subspace of $\mathcal{Y}$ invariant with respect to $\tilde{\mathscr{H}}_h$.\\\\
It is clear that for $u_h\in\mathcal{X}_h^p$ we have
$$
(\tilde{\mathscr{H}}_h u_h,u_h)=\lambda^k |u_h|^2 > 0,\; k=1,\, 2,\ldots.
$$
Let us now consider the restriction $\tilde{\mathscr{H}}_h^p$ of $\tilde{\mathscr{H}}_h$ to $\mathcal{X}_h^p$ which
can be written as
$$
\tilde{\mathscr{H}}_h^p = \sum_{\lambda\in\sigma^p(\tilde{\mathscr{H}}_h)} \lambda E^{\tilde{\mathscr{H}}_h,p}(\lambda),
$$
where the spectral measure $E^{\tilde{\mathscr{H}}_h,p}$ is the point measure, i.e. the orthogonal projection onto $ker(\tilde{\mathscr{H}}_h-\lambda\mathbf{I})$.
With this construction $\tilde{\mathscr{H}}_h^p$ is a positive and self-adjoint operator on $\mathcal{X}_h^p$ with compact inverse $(\tilde{\mathscr{H}}_h^p)^{-1}$.
By Lemma 1, see also Proposition 13.4 in \cite{DAL}, we conclude that there exists a positive and self-adjoint operator $\tilde{\mathscr{H}}^p$ such that, up to a subsequence,
$\tilde{\mathscr{H}}_h^p$
G-converges to $\tilde{\mathscr{H}}^p$, where  $\tilde{\mathscr{H}}^p$ has domain $\mathbf{D}(\tilde{\mathscr{H}}^p) = \mathcal{X}^p$
where
$$
\mathcal{X}^p = \oplus_{k\in\mathbb{N}}\mathscr{N}^k
$$
is a closed subspace of $\mathcal{Y}$ and where
$$
\mathscr{N}^k = \{u\in \mathcal{X};\tilde{\mathscr{H}}^p u=\lambda^k u\}.
$$
Moreover, by Theorem 6, the sequence of $k^{th}$ eigenvalues $\{\lambda_h^k\}$ associated to the sequence $\{\tilde{\mathscr{H}}_h^p\}$
converges to the $k^{th}$ eigenvalue of $\lambda_h^k$ of $\tilde{\mathscr{H}}^p$ and  the corresponding sequence $\{u_h^k\}$
converges to $u^k$ weakly in $\mathcal{X}$.
The limit shifted Dirac operator restricted to $\mathcal{X}^p$ is explicitly given by
$$
\tilde{\mathscr{H}}^p= (\tilde{\mathbf{H}} + V_1M(V_2))|_{\mathcal{X}^p}.
$$
This follows by standard arguments in homogenization theory, see e.g. \cite{BLP}.
\\\\
We continue now to study the asymptotic analysis of the projection to the closed subspace of $\mathcal{Y}$ corresponding to the
positive part $[2,+\infty)$ of the continuous spectrum of $\tilde{\mathscr{H}}_h$.\\\\
We denote by $\mathcal{X}_h^c$ the orthogonal complement in $\mathcal{X}$ to the eigenspace $\mathcal{X}_h^p$. Thus, $\mathcal{X}_h^c$
is the closed subspace invariant with respect to $\tilde{\mathscr{H}}_h$ corresponding to the absolutely continuous spectrum
$\sigma^c(\tilde{\mathscr{H}}_h)=(-\infty,0]\cup[2,\infty)$. We now define the two mutually orthogonal subspaces
$\mathcal{X}_h^{c,+}$ and $\mathcal{X}_h^{c,-}$ with
$$
\mathcal{X}_h^c = \mathcal{X}_h^{c,+}\oplus\mathcal{X}_h^{c,-}
$$
where $\mathcal{X}_h^{c,+}$ corresponds to the positive part $[2,+\infty)$ and $\mathcal{X}_h^{c,-}$ corresponds to the negative part $(-\infty,0]$, respectively.
Next we define the restriction $\tilde{\mathscr{H}}_h^{c,+}$ of $\tilde{\mathscr{H}}_h$ to $\mathcal{X}_h^{c,+}$ which
can be written as
$$
\tilde{\mathscr{H}}_h^{c,+} = \int_{\lambda\in\sigma^{c,+}(\tilde{\mathscr{H}}_h)} \lambda dE^{\tilde{\mathscr{H}}_h,c,+}(\lambda),
$$
where the spectral measure $E^{\tilde{\mathscr{H}}_h,c,+}(\lambda)$ is the continuous spectral measure corresponding to $\tilde{\mathscr{H}}_h^{c,+}$.
By construction $\tilde{\mathscr{H}}_h^{c,+}$ is a positive and self-adjoint operator on $\mathcal{X}_h^{c,+}$. Therefore by Proposition 13.4 in \cite{DAL},
there exists a sequence $\{\tilde{\mathscr{H}}_h^{c,+}\}$ which G-converges to a positive and self-adjoint operator $\tilde{\mathscr{H}}^{c,+}\in\mathcal{X}^{c,+}$.
Moreover, since $\lambda$ is not an eigenvalue, the corresponding sequence $\{E^{\tilde{\mathscr{H}}_h,c,+}(\lambda)\}$ of spectral measures
converges to the spectral measure $E^{\tilde{\mathscr{H}},c,+}(\lambda)$ corresponding to $\tilde{\mathscr{H}}^{c,+}$.\\

Let us consider the evolution equation
\begin{equation}
\label{Evo1}
\left\{\begin{array}{l}
i\frac{\partial}{\partial t}\mathbf{u}_h(t,x) = \tilde{\mathscr{H}}_h^{c,+}\mathbf{u}_h(t,x)\,,\\
\mathbf{u}_h(\cdot,0)=\mathbf{u}_h^0\,.
\end{array} \right.
\end{equation}
By the Stone theorem, there exists a unique solution $\mathbf{u}_h=\mathbf{u}(x,t)$ to (\ref{Evo1}) given by
\begin{equation}\no
\mathbf{u}_h(\cdot,t)=\mathcal{U}_h(t)\mathbf{u}_h^0\, , \;\forall \mathbf{u}_h^0\in \mathcal{X}_h^{c,+}\,,
\end{equation}
where $\mathcal{U}_h(t)=\exp(-i\tilde{\mathscr{H}}_h^{c,+} t)$ is the strongly continuous unitary group of transformations generated by the infinitesimal operator $-i\tilde{\mathscr{H}}_h^{c,+}$ on $\mathcal{Y}$.
By the G-convergence of the sequence $\{\tilde{\mathscr{H}}_h^{c,+}\}$ it follows that the associated sequence $\{\mathcal{U}_h^{c,+}(t)\}$ of unitary groups of transformations
converges to a unitary group of
transformations $\mathcal{U}^{c,+}(t)$ which for every $\mathbf{u}^0\in \mathcal{X}^{c,+}$
defines the solution $\mathbf{u}(\cdot,t)=\mathcal{U}(t)\mathbf{u}^0$ to the limit evolution equation
\begin{equation} \no
\left\{ \begin{array}{l}
i\frac{\partial}{\partial t}\mathbf{u}_(t,x) = \tilde{\mathscr{H}}^{c,+}\mathbf{u}(t,x)\,,\\
\mathbf{u}(\cdot,0)=\mathbf{u}^0\,.
\end{array} \right.
\end{equation}

Finally, by considering the operator $-\tilde{\mathscr{H}}_h^{c,-}$ where $\tilde{\mathscr{H}}_h^{c,-}$ is the restriction to $\mathcal{X}_h^{c,-}$, i.e. the closed subspace corresponding to the negative part $(-\infty,0]$
of the continuous spectrum we can repeat all the arguments from the positive part of the continuous spectrum.

\end{document}